# Sequential tests and estimates after overrunning based on $p$-value combination

## W. J. Hall[*1] and Keyue Ding[†2]

*University of Rochester and Queen's University*

**Abstract:** Often in sequential trials additional data become available after a stopping boundary has been reached. A method of incorporating such information from overrunning is developed, based on the "adding weighted Zs" method of combining $p$-values. This yields a combined $p$-value for the primary test and a median-unbiased estimate and confidence bounds for the parameter under test. When the amount of overrunning information is proportional to the amount available upon terminating the sequential test, exact inference methods are provided; otherwise, approximate methods are given and evaluated. The context is that of observing a Brownian motion with drift, with either linear stopping boundaries in continuous time or discrete-time group-sequential boundaries. The method is compared with other available methods and is exemplified with data from two sequential clinical trials.

## Contents

## 1. Introduction

Suppose a sequential trial is carried out to test a null hypothesis about a real parameter $\delta$. Once the trial is concluded, a non-sequential trial is conducted, with a test of the same hypothesis. The trials are connected in that the amount of

---

*Supported in part by grant R01-HL58751 from the National Heart, Lung and Blood Institute (USA).

†Supported by a grant from the Natural Sciences and Engineering Research Council of Canada.

¹Department of Biostatistics, University of Rochester Medical Center, Rochester, NY 14642-8630, USA, e-mail: hall@bst.rochester.edu

²NCIC Clinical Trials Group, Queen's University, Kingston, ON K7L 3N6, Canada

33

information in the non-sequential trial may depend on data accumulated in the sequential trial. How can the results of the two trials be combined, and a single overall test constructed? The context is that the data, or incremental information, in the non-sequential trial represent "overrunning," from "lagged" data from the sequential trial.

T. W. Anderson [1] considered the problem of incorporating lagged data in an accept-reject rule following a *sequential probability ratio test* and proposed an (approximate) *likelihood ratio test*. In the context of modern-day clinical trials, the problem of how to incorporate data from overrunning was raised and discussed by Whitehead [16, 17], and he gives an admittedly *ad hoc* solution, later named the *deletion method* [14]. This latter paper includes a comparison of the *deletion method* with methods described herein, under certain limited conditions. Another solution is presented in Hall and Liu [4] – actually, an extension of Anderson's likelihood ratio method – along with a discussion of the possible structure of overrunning information in a sequential clinical trial. However, this solution utilizes the *maximum-likelihood ordering* of the sample space, requiring specification of the details of the stopping rule beyond the time a stopping boundary was first reached, in contrast to *stagewise ordering*. In this paper, we focus on procedures that do not require such specification. See these references for further introductory material.

In the context of monitoring a Brownian motion with drift by periodic observations – the context considered herein – Whitehead [16] proposes treating the final analysis that incorporates the overrunning data as if it were a scheduled analysis, but ignoring the analysis that led to stopping, and hence involving a *deletion*. He uses a *stepwise ordering* (as defined in [6], for example) for computing *p*-values and carrying out further inference.

Here we provide another solution, based on the "adding weighted Zs" method of combining *p*-values (Stouffer et al. [15], Mosteller and Bush [10], Liptak [7]); one *p*-value is derived from the sequential experiment (without the overrunning) and the other is based solely on the incremental overrunning data. We recommend weighting the two *p*-values using observed information. This is fully legitimate *only* if (i) the amount of information in the non-sequential trial (overrunning) is proportional to that available at termination of the sequential trial or (ii) the sequential trial was actually nonsequential (and test statistics are normally distributed). For discussion of (i), see [4], Section 2.

Another issue that arises in popular group-sequential trials is that if stopping does not occur until the last scheduled analysis, such an analysis will ordinarily not be done until lagged data are available, in which case a *p*-value will be computed by standard group-sequential methods with a re-scheduled final analysis. (This is consistent with the *deletion method*.) A modification of our method, which combines *p*-values for such trials only when stopping early, is evaluated numerically.

Another application of the combination method could be to a *double sampling* study in which the second sample size depends on the outcome – e.g., on the observed variability – of the first sample. These methods are also appropriate for a *meta-analysis* of two (or more) experiments, whether sequential or not.

Brannath, Posch and Bauer [2] proposed *p*-value combination rules in a different context, namely that of *adaptive group-sequential sampling*. In their setting, allowance is made for the possibility of not carrying out the second stage. (In our context, this would constitute "preventing overrunning.") If the second stage is carried out, the two *p*-values are combined in a way that (i) preserves an overall significance level and (ii) recognizes the stopping rule. As here, the second stage *p*-values may be conditional on results from the first stage. They extend to multiple stages recursively. Numerical integration may be required.

The "adding weighted Zs" combination method is described in Section 2 and extended to an ordered sequence of possibly dependent experiments. In Section 3, this method is applied to sequential clinical trials with overrunning. Special attention is given to the case of a constant amount of overrunning information, the case considered in [14], or to an amount proportional to the available at the end of the sequential trial. It is shown that the latter assumption justifies the use of weights related to the observed (and hence random) amounts of information. Otherwise, the use of such random weights leads to a null distribution of the *p*-value which is only approximately uniform. Still, we recommend this usage so long as the approximation is adequate.

In Section 4 we show how to use the *combination p-value method* to compute estimates and confidence intervals, and in Section 5 provide formulas for evaluating the true confidence coefficients associated with these methods, thus enabling an evaluation of approximations noted above. Some evaluations are summarized in Section 6.

In Sooriyarachchi et al. [14], the issue of reversals in the conclusions after incorporating overrunning, from rejection to acceptance of a null hypothesis or vice versa, was raised. They found, in the cases treated numerically there, that both the *deletion method* and the *combination method* might lead to an uncomfortable level of reversals, with the *deletion method* doing so less frequently. They also noted that both methods (in cases treated) sometimes lead to reduced power. We consider these issues in Section 7 and indicate a modification of the *combination method* that reduces these effects.

The methods are applied to data from the MADIT trial [8] in Section 8 – for both the actual linear-boundary design and for an imagined group-sequential version. Results are compared with those from the *deletion* and *ML-ordering methods*. Results from a second example [9] are briefly summarized.

Some final comments appear in Section 9, including a summary comparison of the alternative methods for incorporating overrunning.

## 2. Combining *p*-values by adding weighted Zs and an extension

We suppose some potential data $X$ are to be available for testing a null hypothesis about a real parameter $\delta$ belonging to an interval $\Delta$. For each $\delta_o \in \Delta$, we consider a test of $\delta = \delta_o$ versus $\delta > \delta_o$, with *p*-value $p(x; \delta_o)$ when $X = x$ is observed. Suppose, for each $\delta_o$, $P \equiv p(X; \delta_o)$ is uniformly distributed on (0,1) when $\delta = \delta_o$ and that, for each $x$, $p(x; \delta)$ is increasing in $\delta$. Then $\mathcal{P} \equiv \{p(\cdot; \delta_o) | \delta_o \in \Delta\}$ defines a *proper family of p-values* for this testing problem.

This is an overly strict definition. We have restricted attention to test functions with continuous distributions, and *stochastic ordering* (increasing or nondecreasing) of *p*-values would allow for differing sample spaces, but the conditions given meet our application. We usually omit the word "proper."

The ordering is needed to avoid possible inconsistencies. Data considered to be "more extreme" than that observed should have higher probability under an alternative hypothesis than under the null. Moreover, it facilitates construction of consistently defined confidence bounds. Simply equate the *p*-value for testing $\delta_o$ to $\gamma$ $(1 - \gamma, \text{resp.})$ and solve for $\delta_o$ to obtain a lower (upper, resp.) confidence bound with confidence coefficient $1 - \gamma$. Choosing $\gamma = 0.5$ yields a median-unbiased estimate.

Suppose $p_1$ and $p_2$ are independent $p$-values for the same null hypothesis, and let $z(u) \equiv \bar{\Phi}^{-1}(u)$ with $\bar{\Phi}(z) = 1 - \Phi(z)$ and $\Phi$ being the standard normal distribution function. Let $w_1$ and $w_2$ be positive numbers for which $w_1^2 + w_2^2 = 1$. Then

$$(2.1) \qquad\qquad p \;\equiv\; \bar{\Phi}\,(w_1\,z(p_1) \;+\; w_2\,z(p_2))$$

is the *adding weighted Zs combined p-value* [7, 10, 12, 15]; also see [13]. It is readily seen that the argument of $\bar{\Phi}$ in (2.1), with $p_i$ replaced by the random variable $P_i$, is distributed as standard normal under the null hypothesis, and hence $P$ in (2.1) is distributed as $U(0,1)$. Moreover, $p$ tends to be small whenever both $p_1$ and $p_2$ are small. More precisely, it is seen to be proper whenever $p_1$ and $p_2$ are proper – the $p_i$'s are increasing in $\delta_o$, so the $z(p_i)$'s are decreasing, as is any positively-weighted linear combination, and hence $p$ is increasing in $\delta_o$.

The interpretation should be clear: $z(p_i)$ is a standardized normal deviate that corresponds to the test statistic on which $p_i$ is based (whether or not $p_i$ was based on a normally distributed statistic), and the argument of $\bar{\Phi}$ in (2.1) represents a (weighted) pooling of normal deviates for the two independent tests, with $p$ the resulting $p$-value.

This *combination method* may be extended to settings where $p_1$ and $p_2$ are derived from overlapping data sets but $p_2$ is a conditional $p$-value for each subset of data on which $p_1$ is based. A possible context is that a second experiment was designed based on the outcome of the first experiment, and a conditional test was used in the second experiment. Formally,

**Proposition 2.1.** *Suppose $p_1 \equiv p_1(x;\delta)$ and $p_2 \equiv p_2(x,y;\delta)$, and $\mathcal{P}_1 \equiv \{p_1(\cdot;\delta_o)|\delta_o \in \Delta\}$ is a family of $p$-values and $\mathcal{P}_2 \equiv \{p_2(x,\cdot;\delta_o)|\delta_o \in \Delta\}$ is, for each $X = x$, a family of conditional $p$-values. Then $\mathcal{P}_2$ is a family of unconditional $p$-values, $p_1(X,\delta_o) \perp p_2(X,Y;\delta_o)$ (independent) for each $\delta_o$, and (2.1) defines a family of $p$-values.*

*Proof.* Since $p_2(X,Y;\delta_o)$ is conditionally $U(0,1)$ for every $X$, it is unconditionally $U(0,1)$, and the needed monotonicity also follows. For each $\delta_o$, the joint distribution function of $(P_1, P_2)$ is

$$\begin{aligned}
\Pr\{P_1 \le u_1,\, P_2 \le u_2\} \;&=\; E\big\{1(p_1(X) \le u_1) \cdot E\big[1(p_2(Y,X) \le u_2 \,\big|\, X\big]\big\} \\
&=\; E\big\{1(p_1(X) \le u_1) \cdot u_2\big\} \;=\; u_1\,u_2\,,
\end{aligned}$$

from which independence follows, and this is sufficient for the claim about $p$. □

Now what about the weights? Ordinarily, they might be related to sample size or information. Specifically, if the $p_i$'s are derived from tests based on means of $n_i$ normally distributed observations (with common variance), then a combined $p$ with $w_i \propto \sqrt{n_i}$ would yield the same $p$ as that from a pooling of the two samples. So far, we have only assumed the weights to be positive constants – depending neither on $\delta_o$ nor on the data. Here are some partial extensions; examples of each appear in the next section.

*The weights may depend on $\delta_o$* without affecting the null distribution of $P$ in (2.1), but the monotonicity in $\delta_o$ may be destroyed except for special choices.

*The weights may be random* (depending on $X$) without affecting the monotonicity in $\delta_o$, but would typically disturb the uniformity of the null distribution of $P$.

It should be emphasized that all $p$-values considered above are for one-sided alternatives. After including overrunning, the usual convention of doubling them for 2-sided alternatives may be appropriate.

Finally, we note that all of this can be directly extended to an ordered set of several *p*-values, each involving new data and conditional on all past data, and combined in the "adding weighted Z's" fashion.

Specifically, let $p_k$ be a *p*-value for the incremental stage-$k$ data, conditional on data from all prior stages. Then define a stage-$k$ combination *p*-value by replacing the argument of $\bar{\Phi}$ in (2.1) by $\sum_{i=1}^{k} w_{k:i} z(p_i)$ with stage-$k$ weights all positive, satisfying $\sum_{i=1}^{k} w_{k:i}^2 = 1$, and $w_{k:i}^2 = w_{k-1:i}^2 \cdot (1 - w_{k:k}^2)$ for $i < k$. Equivalently, (2.1) may be applied recursively, replacing $p_1$ by a combined $p$ from earlier stages with weight $w_1$ for this new $p_1$ and $w_2$ for the incremental data, with $w_1^2 + w_2^2 = 1$.

## 3. Incorporating overrunning by combining *p*-values

We now assume a sequential experiment takes place, resulting in an observation of $(T, X)$, say. We focus on the context of observing a Brownian motion $X(t)$ with drift $\delta$, with a stopping time $T$ and $X \equiv X(T)$ upon stopping, but other contexts may be treated similarly. After stopping, some additional data become available, represented by further observation of the process for $t_o = t_o(T, X)$ units of time. Conditional on $t_o$, a sufficient statistic for the overrunning data is the increment $Y$ observed during the overrunning time increment $t_o$. In other words, a sequential experiment is followed by a non-sequential one, with sample size (observation time) depending on the outcome of the sequential trial. There may be additional randomness in $t_o$; it is sufficient to let $t_o(t, x)$ be the conditional expectation of overrunning information, given $(T, X) = (t, x)$. See ([4], Section 2) for discussion supporting $t_o$ being a constant, $\propto \sqrt{t}$, or $\propto t$ as possible approximations to reality.

Upon reaching a stopping boundary, a *p*-value $p_1$ for a null hypothesis about the drift parameter is defined: $\delta = \delta_o$ versus $\delta > \delta_o$. And at the end of overrunning, a conditional *p*-value $p_2$ is simply $\bar{\Phi}\big((y - \delta_o t_o)/\sqrt{t_o}\big)$, given $t_o = t_o(t, x)$. A combination *p*-value is therefore given by (2.1). (Here, $(T, X)$ plays the role of $X$ in Section 2.) Hence,

**Corollary 3.1.** *Suppose $w_1$ and $w_2$ are positive constants for which $w_1^2 + w_2^2 = 1$. Then*

$$(3.1) \qquad p(t, x, y; \delta_o) \;\equiv\; \bar{\Phi}\left(w_1\, z(p_1(t, x; \delta_o))\; +\; w_2(y - \delta_o t_o)/\sqrt{t_o}\,\right)\Big|_{t_o = t_o(t, x)},\; \delta_o \in \Delta\,,$$

*defines a family of p-values.*

But how should the weights be chosen? It is tempting to choose them to be proportional to the square-root of *information* in the respective parts of the experiment. Then each summand in (3.1) would have variance or conditional variance equal to the information in that part of the experiment. Using expected information, $w_1^2 = E_{\delta_o}(T)/[E_{\delta_o}(T) + E_{\delta_o} t_o(T, X)]$ and $w_2^2 = 1 - w_1^2$. But, as noted in Section 2, this would not typically preserve the needed monotonicity of $p(\delta_o)$ in (3.1). Moreover, knowledge of the functional form of the dependence of $t_o$ on $(t, x)$ would be needed. If $t_o$ were constant, this would yield

$$p(t, x, y; \delta_o) \;=\; \bar{\Phi}\left(\frac{[E_{\delta_o}(T)]^{1/2}\, z\big(p_1(\delta_o)\big)\; +\; y - \delta_o t_o}{[E_{\delta_o}(T)\; +\; t_o]^{1/2}}\right).$$

This could be used as a *p*-value for a single null hypothesis, but it would not be suitable for construction of confidence bounds, unless $E_{\delta_o}(T)$ was replaced by $E_{\delta_o'}(T)$ for a fixed $\delta_o'$. Because of these limitations, we abandon this approach.

Suppose instead we use the square-roots of *observed information*, namely $\sqrt{t}$ and $\sqrt{t_o}$, yielding

$$(3.2) \qquad p(t, x, y; \delta_o) \; = \; \bar{\Phi}\left( \frac{t^{1/2}\, z\big(p_1(t, x; \delta_o)\big) \, + \, y - \delta_o t_o(t, x)}{[t \, + \, t_o(t, x)]^{1/2}} \right).$$

Monotonicity in $\delta_o$ (for each $(t, x, y)$) is maintained, but the uniformity of the null distribution would appear to be in doubt. However, to compute $p$, no knowledge of the dependency structure of $t_o$ is required, only its observed value.

We now consider the special case of (3.1) and (3.2) with $t_o \propto t$, say $t_o = c\,t$. Since $w_1^2 = t/(t + ct) = 1/(1 + c)$ and $w_2^2 = c/(1 + c)$, this yields constant weights; and $c$ is known once $T = t$ and $t_o$ are observed. Hence, the use of observed information in this case is justified.

**Corollary 3.2.** *If, for some constant $c$, $t_o(x, t) = ct$ for all $(x, t)$, then (3.2) defines a family of p-values.*

For the group-sequential case with up to $K$ analyses and stagewise ordering, we modify the combination $p$-value (3.2): For testing $\delta_o$, with $t_{ok} \equiv t_o(k)$, the modified $p$-value is defined as

$$(3.3)$$
$$p^*(t_k, x, y; \delta_o) = \begin{cases} \bar{\Phi}\left([t_k^{1/2}\, z(p_1(t_k, x; \delta_o)) + y - \delta_o t_{ok}]/(t_k + t_{ok})^{1/2}\right) & \text{if } k < K \\ p_1(t_K + t_{oK}, x + y; \delta_o) & \text{if } k = K \end{cases}$$

where $p_1(t_k, x; \delta_o)$ is the group-sequential stagewise $p$-value for testing $\delta = \delta_o$ versus larger values when the analyses are scheduled at $t_1, \ldots, t_{K-1}, t_K^o \equiv t_K + t_{oK}$ with early-stopping sets $\mathcal{S}_k$ $(k < K)$ (each the complement of an interval). This matches the *deletion method* when stopping has not occurred early.

For a group-sequential ML-ordering, the *ML-ordering method* [4] may be more suitable.

We show in Sections 5 and 6 that use of $p$ in (3.2) or (3.3), for several choices of the dependency of $t_o$ on $(t, x)$ and two popular sequential designs, for constructing confidence bounds and intervals may yield adequately accurate confidence coefficients. This leads us to recommend the use of (3.2) or (3.3) as if it were a bona fide combination $p$-value, if the design chosen and the likely form of dependency are similar to those considered in Section 6.

One last variation permits further adjustment of the weighting: Use weights with squares proportional to $T$ and $\rho\, t_o(T)$ for a specified weighting factor $\rho > 0$. For motivation, see Section 7.

## 4. Computing *p*-values and confidence bounds

Here we act as if $t_o \propto t$, and discuss the use of (3.2) and (3.3) for obtaining $p$-values and, by inversion, confidence bounds and intervals.

For any particular null value $\delta_o$, the *combined p-value* $p(\delta_o)$ may be computed from (3.2) or (3.3) with $t$ (or $t_k$), $x$, $y$ and $t_o$ the observed values, and using software that enables computation of $p_1(\delta_o)$. For general linear boundaries, such software is available from the authors (based on formulas in [3]), and the PEST software [11] provides such output for a limited selection of linear boundaries and group-sequential modifications of them. For group-sequential boundaries with stagewsie ordering, a program – built around software for $p_1(\delta_o)$ from Jennison [5] – is available from the authors.

To obtain an upper confidence bound with confidence coefficient $\gamma$, we need to solve $p(\delta) = \gamma$ for $\delta = \hat{\delta}_U$, or equivalently, solve $z(p(\delta)) = z(\gamma)$. A little algebra leads to the equivalent problem – except in the group-sequential case with $t = t_K$ – of solving $\delta - h(\delta) - [y - \sqrt{t^o}\, z(\gamma)]/t_o = 0$ where $t^o \equiv t + t_o$ and $h(\delta) \equiv \sqrt{t}\, z(p_1(\delta))/t_o$. Starting from a trial solution $\delta^o$, and computing $h(\delta^o)$ and $h(\delta^o + \epsilon)$ for some small $\epsilon$, an improved solution is

$$\delta \;\equiv\; \delta^o \;-\; \frac{\delta^o - h(\delta^o) - [y - \sqrt{t^o}\, z(\gamma)]/t_o}{1 \;+\; [h(\delta^o) - h(\delta^o + \epsilon)]/\epsilon}\,.$$

We find that two or three iterations provide good accuracy. (When $t = t_K$, we only need solve $p_1(t_K^o, x + y; \delta) = \gamma$.) Alternatively, a trial-and-error approach works quite satisfactorily.

## 5. True confidence coefficients

We now evaluate the true confidence coefficient for a confidence bound or interval determined by using (3.2), whether or not $t_o$ is proportional to $T$. Let $\hat{\delta}_\gamma$ be an upper confidence bound determined by the method of the previous section for a nominal confidence coefficient $\gamma$. The question is: what is the true confidence coefficient? We need to evaluate, for given $\gamma$ and $\delta$, $q_\gamma(\delta) \;\equiv\; P_\delta(\delta < \hat{\delta}_\gamma)$ and determine $q_\gamma^o \equiv \inf_\delta q_\gamma(\delta)$.

As noted in Section 4, $\hat{\delta}_\gamma$ is the solution to $\delta - h(\delta) = [y - \sqrt{t^o}z(\gamma)]/t_o \equiv g(y, t^o, t_o, \gamma)$. Since $h$ is decreasing in $\delta$, the left side is increasing in $\delta$, and hence $\delta < \hat{\delta}_\gamma$ iff $\delta - h(\delta) < g(y, t^o, t_o, \gamma)$. This latter event is equal to the event $(y - \delta t_o)/\sqrt{t_o} > [-\sqrt{t}\, z(p_1(t, x; \delta)) + \sqrt{t^o}\, z(\gamma)]/\sqrt{t_o}$. Therefore, conditioning on $(T, X)$ and hence on $T_o$, we have

$$(5.1) \qquad q_\gamma(\delta) = E_\delta P_\delta(\delta < \hat{\delta}|T, X) = E_\delta \Phi\left(\frac{T^{1/2}\, z(p_1(T, X; \delta)) - T^{o\,1/2}\, z(\gamma)}{T_o^{1/2}}\right).$$

If this combination $p$-value were bona fide – that is, if $T_o \propto T$ – the result would be $\gamma$ identically in $\delta$.

The true confidence coefficient for an (equal-tail) confidence interval based on (3.2) may be obtained similarly. For an interval with nominal confidence coefficient $\gamma$, the true confidence coefficient is $Q_\gamma^o \equiv \inf_\delta Q_\gamma(\delta)$ where

$$(5.2) \qquad\qquad Q_\gamma(\delta) \;\equiv\; q_{(1-\gamma)/2}(\delta) - q_{(1+\gamma)/2}(\delta)\,.$$

For the group-sequential modification (3.3), (5.1) needs to be modified when $T = T_K$.

## 6. Computational support for approximations

Here we report on some numerical evaluations of the validity of using (3.2) or (3.3) when $t_o$ is not proportional to the observed stopping time $t$, and the validity of using the combination method only when stopping after an interim analysis. For various special cases and many values of $\delta$, we computed (5.1) for $\gamma = 0.5$ and (5.2) for $\gamma = 0.9$ and 0.95 to see how close they are to the respective nominal values of 0.5, 0.9 and 0.95. We summarize some of the findings here.

*A linear-boundary design:* Consider triangular boundaries for testing $\delta = 0$ versus $\delta = 1$ with intercepts $\pm 5.99$, slopes 0.75 and 0.25, and apex at $t = 23.97$. This design has both error probabilities 0.025. The design may be adapted for testing $\delta = 0$ versus $\neq 0$ (as prescribed by the PEST software). The resulting one-sided rejection region is the upper boundary for which the power at $\delta_1 \equiv 0.8233$ is 0.9. The expected stopping time is 7.776 at $\delta = 0$ (or 1) and 9.382 at $\delta_1$, and has its maximum of 11.217 at $\delta = 0.5$.

We considered $t_o \propto T$, $t_o$ constant and $t_o \propto \sqrt{T}$. For the first case, we simply verified the accuracy of our computer program, finding that the distribution of the $p$-value was exactly uniform, and that the true confidence coefficients matched the nominal ones exactly.

For the constant case, we considered $t_o = c\,E_{\delta_1}(T)$ with $c$ ranging from 0.1 to 0.5. Here are selected results:

| $c = 0.1$ | $c = 0.5$ |
|---|---|
| $0.487 < q_{.5} < 0.513$ | $0.471 < q_{.5} < 0.529$ |
| $0.900 < Q_{.9} < 0.908$ | $0.899 < Q_{.9} < 0.917$ |
| $0.950 < Q_{.95} < 0.955$ | $0.949 < Q_{.95} < 0.960$ |

For $c = 0.1$, the true confidence coefficients $Q_\gamma^o$ for nominal 90% and 95% confidence intervals are therefore correct (to 3 decimal places), and only slightly below the nominal values for $c = 0.5$. However, the median-unbiased estimate may have a few percentage points of median-bias, depending on the true $\delta$. We also found that $q_{.5} < 0.5$ for $\delta > 0.5$ and vice versa. Computations for $c$-values between 0.1 and 0.5 yielded bounds between the respective ones in the display above. Results for $t_o \propto \sqrt{T}$ were uniformly better than those for $t_o$ constant.

*An O'Brien–Fleming group-sequential design:* Consider an O'Brien–Fleming two-sided design for testing $\delta = 0$ with significance level 0.05 and power 0.9 at $\delta = \pm 1$, with a maximum of 5 analyses. We assume equally spaced interim analyses, at 0.2, 0.4, 0.6, 0.8 times $t_5 \equiv 10.781$, with boundary values of $\pm 6.6988$ (obtained from [6]). Again, we considered $t_o \propto T$ for the unmodified combination $p$-value to confirm the accuracy of our programs. For the modified $p$, we considered $t_o$ constant, namely $= c\,t_5$, and $t_o = c\,t_k$; in each case, $c$ ranged from 0.02 to 0.1.

Here are some of the results:

|  | $t_o = c\,t_5$ | | $t_o = c\,t_k$ | |
|---|---|---|---|---|
|  | $c = 0.02$ | $c = 0.1$ | $c = 0.02$ | $c = 0.1$ |
| $q_{.5}^o$ | 0.475 | 0.445 | 0.478 | 0.451 |
| $Q_{.9}^o$ | 0.894 | 0.887 | 0.894 | 0.888 |
| $Q_{.95}^o$ | 0.947 | 0.943 | 0.947 | 0.943 |

Again, although $q_{.5}^o$ may be as small as 0.44 (and by symmetry $0.44 < q_{.5}(\delta) < 0.56$), we found that $q_{.5}$ was usually within $\pm 0.01$ of 0.5. Indeed, this occurred for all but 1%, 7%, 1% and 4%, respectively (reading from left to right in the display above) of the range of $\delta$-values within $\pm 2.5$.

## 7. Reversals and power

Sooriyarachchi et al. [14] raised concern about the frequency of reversals of acceptance and rejection conclusions after inclusion of overrunning information, but stressed their desire not to ignore such information. In simulation studies of the *deletion* and *combined p-value methods*, with constant amounts of lagged data (independent of the results at the time of stopping), they found levels of reversals that they considered worrisome, especially for the *combination method* – perhaps 3 or 4 percent. However, in popular group-sequential designs such as O'Brien–Fleming, reversals were rare and only defined when the trial stopped early, as an analysis at a final scheduled time would ordinarily await lagged data before execution.

Of more concern to us, is their finding that both methods may lead to reductions in power. Intuitively, when a rejection occurs "early", overrunning can reverse it but the chances of compensating with reversals in the other direction may be minimal.

With constant overrunning information, our computations (not reported here) confirm theirs, but we find reversals to be somewhat less frequent when overrunning information increases with stopping times, and losses in power are then rarer.

A possible compromise method is as follows: down-weight the overrunning *p*-value in the combination formula. By introducing a factor $\rho$ (see end of Section 3), it is possible to maintain power and depress the frequency of reversals but still not ignore the lagged data completely. However, computations show that some situations will require extensive down-weighting (small $\rho$). Choice of a suitable $\rho$ will require computational trial-and-error, with assumptions about overrunning needed. For this purpose, we provide the following formulas.

When the true drift is $\delta$ (and stopping is in continuous time), the probability of rejection upon stopping followed by acceptance after inclusion of overrunning, when $t_o \propto t$, is

(7.1)
$$P_\delta(R \to A) = \int_0^{t_{max}} \Phi\Big\{\big[(1+\rho c)/(\rho c)\big]^{1/2} z_\alpha - \big[1/(\rho c)\big]^{1/2} z_1(U,t) - \delta(ct)^{1/2}\Big\} dP_\delta^U(t)$$

with $z_1(U,t)$ being the standard normal deviate for which the right-hand-side tail area beyond it is $P_0^U(t)$ (the *p*-value when the upper boundary $U$ is crossed at time $t$), $P_\delta^U(t)$ being the probability of crossing the upper boundary before the lower one prior to $t$, and $\alpha$ being the one-sided significance level for testing $\delta = 0$. For group-sequential tests, the integrator in (7.1) is $dP_\delta^U(x,t)$, indicating a need to integrate over $x$-values where $t = t_k$ and the upper boundary has been reached, but $t$ may be restricted to $\{t_k | k < K\}$ since reversals at a final analysis have no role.

Similarly, $P_\delta(A \to R)$ is given by (7.1) with $U$ replaced by $L$ (for lower boundary) and $\Phi$ replaced by $\bar{\Phi}$. Finally, the power after inclusion of overrunning, when the power of the original design is $pow(\delta)$, is

$$ovpow(\delta) = pow(\delta) - P_\delta(R \to A) + P_\delta(A \to R).$$

(Software is available from the authors.)

## 8. An example: the MADIT study

MADIT (Multicenter Automatic Defibrillator Implantation Trial [8]) was a randomized clinical trial conducted to evaluate the effectiveness of an implanted defibrillator compared with conventional drug therapy to reduce mortality associated with ventricular arrhythmias. Monitoring was based on the logrank statistic plotted against its estimated variance [17]. This behaves like a Brownian motion with drift $\delta = -\log(HR)$ where HR is the hazard ratio of the treatment-to-control arms (assuming proportional hazards). The essential features were reviewed in [4] and are summarized here.

A triangular design was used that assures a two-sided significance level of 5% and a power of 90% at a hazard ratio of 0.537 (drift = 0.6218). Monitoring was carried out weekly over the five years of the trial, thereby yielding nearly-continuous observation of the logrank process. The stopping boundaries were $u_t = 7.935 + 0.189t$ and $l_t = -7.935 + 0.566t$, with the early part of the lower boundary ($l_t$) a rejection region for superiority of the control arm.

Interpolating, the upper boundary was reached at $t = 12.145$ with $x = 10.230$, later corrected to $t = 12.037$ and $x = 10.210$. The incremental coordinates for overrunning were $t_o = 1.240$ and $y = 2.957$, showing an upturn in the sample path after reaching the boundary.

Respective $p$-values and estimates of the drift and of the HR are presented below, contrasting results of analyses without and with the use of the overrunning data. Values in square brackets are those reported in [4] for the *ML-ordering method*, assuming $t_o \propto \sqrt{t}$; with linear boundaries and no overrunning, stepwise ordering and ML-ordering are identical.

| *overrunning* | *2-sided p* | *med-unb-est* | *95% confidence interval* |
|---|---|---|---|
| Inference about the drift $\delta$ | | | |
| *without* | 0.0084 | 0.786 | (0.204, 1.361) |
| *with* | 0.0009 [0.0029] | 0.938 [0.939] | (0.388, 1.484) [(0.329, 1.543)] |
| Inference about the hazard ratio $HR = \exp(-\delta)$ | | | |
| *without* | 0.0084 | 0.456 | (0.256, 0.815$^+$) |
| *with* | 0.0009 [0.0029] | 0.391 [0.391] | (0.227, 0.678) [(0.214, 0.720)] |

Both methods reflect the upturn in the sample path during overrunning as the "with" $p$-values are smaller and the estimates farther from the null values. But the combination method gives the smaller $p$-value and narrower confidence intervals; this may reflect the different orderings being used by the two methods.

Values reported in [8] were based on Whitehead's *deletion method*; they are identical to the "without overrunning" values in the display above, as the *deletion method* essentially ignores overrunning when the path continues in a similar direction and there is near-continuous monitoring. (It was this observation that inspired the development of alternative methods for incorporating overrunning.) In such settings, the deletion-method $p$-value cannot be smaller than when computed upon first hitting an upper boundary, irrespective of the nature of the overrunning data. (For, when reaching the upper boundary at time $t$, with the prior analysis a short time earlier, at time $t^-$ say, and then overrunning to $x^o$ at a later time $t^o$, the deletion one-sided $p$ is the null probability of $\{T \leq t^- \text{ and } X(T) \geq u_T\} \cup \{T \geq t \text{ and } X(t^o) \geq x^o\}$. These two events are disjoint, and the former is virtually the extremal set without overrunning.)

We now consider a group-sequential variation on MADIT as described in [4]. We pretended that an O'Brien–Fleming 5-analysis design was used for testing $\delta = 0$ versus $\neq 0$ with power 80% at a HR of 0.537. It would have stopped at the third interim analysis with the results obtained upon hitting the boundary in MADIT. Results of analyses are reported below; for comparison, values in square brackets are those reported in Table 2 of [4] using the *ML-ordering method* and assuming $t_o \propto \sqrt{t}$. Values for the *deletion method* – which treats the analysis after overrunning as a replacement for the third scheduled analysis – are also given.

In each case – i.e., without or with overrunning – results from the group-sequential combination method indicate a more significant departure from the null value of $HR = 1$ than do those by the group-sequential ML-ordering method. At least for the "without' results, this is attributable to the different orderings used. This time the *deletion method* gives results similar to those from *ML-ordering.* (These results are not directly comparable to those in the previous table since the pretended group-sequential design has reduced power.)

| *overrunning* | *2-sided p* | *med-unb-est of HR* | *95% confidence interval* |
|---|---|---|---|
| Group-sequential inference about the hazard ratio | | | |
| *without* | 0.0039 [0.0041] | 0.431 [0.468] | (0.244, 0.762) [(0.295, 0.775)] |
| *with* | 0.0004 [0.0011] | 0.373 [0.384] | (0.217, 0.641) [(0.221, 0.672)] |
| Deletion method | 0.0014 | 0.384 | (0.221, 0.680) |

Here is a brief summary of results from a second defibrillator trial, MADIT-II [9], in which the *combination method* was pre-specified. The design was again triangular, with a 5% 2-sided significance level and power 95% at a HR of 0.627: $u_t = 11.77 + 0.1273\,t$ and $l_t = -11.77 + 0.3819\,t$. This time $(t, x, t_o, y) = (45.415, 17.551, 0.483, 1.441)$. The results were:

| *overrunning* | *2-sided p* | *med-unb-est of HR* | *95% confidence interval* |
|---|---|---|---|
| Inference about the hazard ratio in MADIT-II | | | |
| *without* | 0.028 | 0.708 | (0.525, 0.962) |
| *with* | 0.016  [0.023] | 0.688 [0.689] | (0.511, 0.932) [(0.504, 0.948)] |

Again, the *deletion method* of incorporating overrunning would have agreed with the "without" analysis, and results from the *ML-ordering method* (in square brackets) are mainly intermediate.

## 9. Final remarks

Proposition 1 applies to other methods of combining *p*-values, such as Fisher's summing of $-\log(1 - p_i)$. (For a description of such methods, see [12] or [13].) We chose the "adding Zs" method for two reasons: (i) It lends itself naturally to weights – it would be unreasonable to give equal weights to a long trial and a small amount of overrunning – and (ii) it reduces to standard normal-theory methods when the sequential component is replaced by a non-sequential one – equivalently, if a naive analysis is done after stopping rather than one recognizing the stopping rule.

Here are some of the pros and cons of various methods for incorporating overrunning:

(a) *Deletion method:* Not suitable for near-continuous monitoring. Ignores the fact that, at the boundary-hitting stage, the monitoring statistic was in a stopping region but is the natural approach in a group-sequential trial when early stopping has not occurred. Simple to use. Results in approximate *p*-values and final inference. Limited computations show that a loss in power may occur.

(b) *Combination p-value method:* Makes direct use of the analysis that led to stopping. Approximate except when $t_o \propto T$, and even then for common group-sequential designs. Uses stage-wise ordering, and hence free of any direct dependence on future stopping boundaries. Needs no formal assumption about the form of $t_o$.

Computations show a loss in power may occur. May be modified to reduce the chance of reversal after overrunning and loss in power.

(c) *ML-ordering method:* Based on a minimal sufficient statistic, and hence ignores which boundary was first reached and when. Exact, up to needed assumptions about overrunning information (and Brownian motion approximation). Requires an assumption about the form of $t_o(t)$, but not very sensitive to it in the practical cases examined. Uses ML-ordering and hence depends on stopping boundaries beyond those when boundaries were first reached.

Sooriyarachchi et al. [14] conclude that (a) is preferable although they only considered constant amounts of overrunning whereas our focus has been on settings where the amount of overrunning information is likely to increase with increased stopping times. They highly stress the possibilities of reversals, but such possibilities cannot be avoided once one agrees to utilize lagged data. The chances can be reduced within the *combination method* by reducing the weight given to lagged data, but this would need to be considered in advance of the trial.

We recommend (c) in settings where the design is likely to be followed closely. A numerical study of reversals and power with the *ML-ordering method* will be presented elsewhere. Otherwise, we think the *combination p-value method*, possibly with a down-weighting of overrunning information, is a competitor worthy of consideration, especially when overrunning increases with increasing stopping times.

We encourage investigation of the *combination method* in other settings, including meta-analyses and double sampling.

**Acknowledgments.** The first MADIT trial stimulated the research reported here. We thank Arthur Moss, MD, and Boston Scientific Corporation (formerly CPI/Guidant), respectively, for leadership and sponsorship of this trial. We are also grateful to John Whitehead for helpful discussion and to Michael McDermott for some references to the *p*-value literature.